\documentclass[12pt]{article}
\usepackage{amsmath}
\usepackage{amssymb,latexsym}

\begin{document}
\title{On Partitions of Finite vector Spaces}
\author{Antonino Giorgio Spera\thanks{Research supported by University of Palermo (FRSA).}\\
Universita' degli studi di Palermo\\
Dipartimento di Matematica ed Applicazioni\\
Via Archirafi 34,  I\-90123 Palermo, Italy\\
\texttt{spera@math.unipa.it}}
\date{   }
\maketitle
\begin{abstract}
\noindent In this note, we give a new necessary condition for the
existence of non-trivial partitions of a finite vector space.
Precisely, we prove that the number of the subspaces of minimum
dimension $t$ of a non-trivial partition of $V_{n}{(q)}$ is greater
than or equal to $q+t$. Moreover,  we give some extensions of a well
known Beutelspacher-Heden's result on existence of $T$-partitions.
\bigskip


\noindent Math. Sub. Cl. (2000): Primary 15A03, Secondary  05B30 and
20D60.

\end{abstract}

\section{Introduction} \label{S:intro}

\bigskip

\noindent A partition ${\bf P}$ of the $n$-dimensional vector space
$V_{n}{(q)}$ over the finite field with $q$ elements, is a set of
non-zero subspaces (components) of $V_{n}{(q)}$ such that each
non-zero element of $V_{n}{(q)}$ is contained in exactly one element
of ${\bf P}$.

The interest on the problems of existence, enumeration, and
classification of partitions, arises in connection with the
construction of codes and combinatorial designs. In fact, if ${\bf
P}=\{V_{1},V_{2},.....,V_{r} \}$ is a non-trivial partition of
$V_{n}{(q)}$, then the subspace $W$ of the vector space $V_{1}\times
V_{2} \times ..... \times V_{r}$, which is defined by
$$W:=\{ (y_{1},y_{2},.....,y_{r}) \in V_{1}\times V_{2}
\times ..... \times V_{r} \hspace{0.05in} | \hspace{0.05in}
\sum_{i=1}^{r}y_{i}= \{ 0 \}\},$$
 is a perfect mixed linear code
(see [HS] and [Li]). Moreover, if $\mathbb{B}$ is the set of
subspaces in ${\bf P}$ together with all their cosets, then
$V_{n}{(q)}$ and $\mathbb{B}$ are, respectively, the point set and
the block set of an uniformly resolvable design which admits a
translation group isomorphic to $V_{n}{(q)}$ (see  [DR] and
[ESSSV2]). Furthermore, combinatorial designs can be associated with
certain more general "partitions" (see for instance [Sc1], [Sc2] and
[Sp]).

 If a partition consists of $x_{i}$ components of
dimension $n_{i}$ for each $i=1,2,...,k$, then the non-negative
integers $x_{1},x_{2},.....,x_{k}$ are a solution of the Diophantine
equation $ \sum_{i=1}^{k}(q^{n_{i}}-1)X_{i}= q^{n}-1$. A fundamental
problem about partitions of $V_{n}{(q)}$ is  to give necessary and
sufficient conditions on non-negative solutions of the above
equation, in order that they correspond to  partitions of
$V_{n}{(q)}$.

\noindent In [ESSSV2] the authors gave such conditions in the case
where $q=2$, $k=2$, $n_{1}=2$ and $n_{2}=3$ (see section 2). In
order that, they proved the following two necessary conditions for
existence of partitions:

\noindent  a)  In a non-trivial partition of $V_{n}{(q)}$ the number
of components of minimum dimension must be greater than or equal to
$2$.

\noindent b) If the components of minimum dimension of a partition
of $V_{n}{(2)}$ have dimension 1,  then their number is greater than
or equal to $3$ (see [ESSSV1]).

Another interesting problem is related to existence results on
$T$-partitions, where, if $T$ is a set of positive integers, a
$T$-partition is a partition ${\bf P}$ of $V_{n}{(q)}$ such that $\{
dim V'
 \hspace{0.05in} | \hspace{0.05in} V' \in {\bf P} \} = T$. Clearly
 the existence of a $T$-partition of $V_{n}{(q)}$ implies the
 existence of a positive solution of the equation \break $ \sum_{i=1}^{k}(q^{n_{i}}-1)X_{i}=
 q^{n}-1$ in the case where  $T=\{n_{1},n_{2},......,n_{k}\}$. A. Beutelspacher and O. Heden
proved (see  [Be] and [He]) the existence of a $T$-partition in the
case where $\min T \geq 2$ and $\max T = \frac{n}{2}$.

In this paper, in section 2 we will recall some definitions and some
known results about the existence of partitions of a finite vector
space. In section 3, we will provide a more general necessary
condition on the minimum dimension components  of a partition.
Precisely, we will shown that the number of components of minimum
dimension $t$, of any non-trivial partition of $V_{n}{(q)}$, is
greater than or equal to $\alpha q+t$ where $\alpha$ is a positive
integer. Finally, in section 4, we will extend  the mentioned
Beutelspacher-Heden's existence result for $T$-partition  of
$V_{n}{(q)}$ in the case where the minimum dimension of the
components is $1$ and in some other cases where the maximum
dimension of the components is consistent with $n$.

 \vspace{1in}

\section{Definitions and first results} \label{S:P*}

\bigskip

In this section we recall some basic properties of partitions of
finite vector spaces and some results about the existence of certain
classes of partitions.

 Let $n$ be a positive integer $(n>1)$, $q$ be a
prime power, $\mathbb{F}_{q}$   be the finite field of order $q$ and
$V_{n}{(q)}$ be the $n$-dimensional vector space  over
$\mathbb{F}_{q}$. A set ${\bf P}=\{V_{1},V_{2},......,V_{r}\} $ of
non-zero subspaces of $V_{n}{(q)}$ is a partition of $V_{n}{(q)}$ if
and only if
 $ \cup_{i=1}^{r}V_{i}=V_{n}{(q)}$ and $V_{i} \cap V_{j}=\{ 0 \}$
for every $i,j \in \{ 1,2,...,r \}$ and $i \neq j$. We call ${\bf
P}$ non-trivial in the case where $r \geq 2$. The elements of ${\bf
P}$ are said to be the components of the partition. Let $T$ be a set
of positive integers and ${\bf P}$ a set of disjoint  non-trivial
subspaces (which is not necessarily a partition). ${\bf P}$ is said
to be a $T$-set of subspaces (a $T$-partition if ${\bf P}$ is a
partition) if the map $dim:{\bf P}\rightarrow T$ is surjective,
where $dim(V_{i})$ is the dimension of the subspace $V_{i}$ for each
$V_{i}\in{\bf P} $. Of course if $T=\{n_{1},n_{2},......,n_{k}\}$,
then $1\leq n_{i}\leq n$ for every $i=1,2,......,k$.

\noindent A. Beutelspacher and O. Heden proved (see [Be] and [He])
the following well known existence result.

\bigskip
\noindent{\bf 2.1  Theorem}.  {\it Let
$T=\{n_{1},n_{2},......,n_{k}\}$ be a set of positive integers with
$ n_{1}< n_{2}<......< n_{k}$.  If $n_{1}\geq 2$, then there exists
a $T$-partition of $V_{2n_{k}}{(q)}$.}

\bigskip

\noindent The above  theorem was proved by  Beutelspacher  [Be] in
the case where $n_{1}=2$ and by  Heden [He] for  $n_{1}>2$.

 Now, let $k$ be a positive integer and $x_{1},x_{2},......,x_{k}$
 be non-negative integers.
If ${\bf P}$ is a partition of $V_{n}{(q)}$ which contains $x_{i}$
components of dimension $n_{i}$ for each $i=1,2,...,k$, then we say
that ${\bf P}$ is of type
$$[(x_{1},n_{1}),(x_{2},n_{2}),......,(x_{k},n_{k})]$$

\noindent or that ${\bf P}$ is an
$[(x_{1},n_{1}),(x_{2},n_{2}),......,(x_{k},n_{k})]$-partition of
$V_{n}{(q)}$ (see [ESSSV2]). Note that it is possible to have $
x_{i}= 0 $ for some $1 \leq i \leq k$. Clearly, in such a case,
there are no components of dimension $n_{i}$. However, as we will
see, such notation will be useful when we associate partitions to
non-negative solutions of some Diophantine equation.
 Later on, for a partition of type
$[(x_{1},n_{1}),(x_{2},n_{2}),......,(x_{k},n_{k})]$ we will always
suppose $ 1\leq n_{1}< n_{2}<......< n_{k}\leq n$.

Let  ${\bf P}$ be an
$[(x_{1},n_{1}),(x_{2},n_{2}),......,(x_{k},n_{k})]$-partition of
$V_{n}{(q)}$. Then it is easy  to show the following necessary
condition.

\bigskip

1) \hspace{0.05in} $(x_{1},x_{2},......,x_{k})$ is a non-negative
solution of the Diophantine equation
$$ \sum_{i=1}^{k}(q^{n_{i}}-1)X_{i}= q^{n}-1. \eqno(1)$$

\noindent Moreover, if $V_{i}$ and $V_{j}$ are two distinct
components of ${\bf P}$, then $dim(V_{i}+V_{j})=dimV_{i}+dimV_{j}$
since $V_{i}\cap V_{j}=0$. Hence, the following necessary conditions
are obtained:

\bigskip

2)\hspace{0.05in} If $ i \neq j$ and $x_{i} \neq 0 \neq x_{j}$, then
$ n_{i}+n_{j} \leq n$.

    \hspace{0.06in} \hspace{0.06in} If $2n_{i}> n$, then $x_{i} \leq 1$.

\bigskip

Note that, by Theorem 2.1, the  equation $(1)$ always admits a
non-negative solution when $n=2n_{k}$ and $n_{1} \geq 2$.

Furthermore, if $W$ is a subspace of $V_{n}{(q)}$ of dimension $s$,
then ${\bf P'}=\{V_{i}\cap W \hspace{0.05in} | \hspace{0.05in} V_{i}
\in {\bf P} \hspace{0.05in} and \hspace{0.05in} V_{i}\cap W \neq 0
\}$ is a partition of $W$ and so, if ${\bf P'}$ is of type
$[(x'_{1},n'_{1}),(x'_{2},n'_{2}),......,(x'_{k'},n'_{k'})]$, we
also have

\bigskip

3) \hspace{0.05in} the equation $\sum_{i=1}^{k'}(q^{n'_{i}}-1)X_{i}=
q^{s}-1$
 admits the non-negative

  \hspace{0.05in} \hspace{0.05in} \hspace{0.05in} solution $ (x'_{1}, x'_{2},.....,x'_{k'}) $.

\bigskip

\noindent In [Bu] it was shown that for $s=n-1$ the property  3)  is
a necessary condition which does not follow from conditions 1) and
2).

Now we recall some existence results. The following two theorems can
be found in [Bu], but Theorem 2.2 was known before. In fact, part i)
of it is a well known result on $d$-spreads and part ii) has been
also proved previously  by Beutelspacher in [Be].

\bigskip
\noindent{\bf 2.2  Theorem}.  {\it Let $d$ and $n$ be positive
integers.

\noindent i) If $d$ divides $n$, then there exists  a partition of
$V_{n}{(q)}$ of type $[(\frac{q^{n}-1}{q^{d}-1},d)]$.

\noindent ii) If $d< \frac{n}{2}$, then there exists a partition of
$V_{n}{(q)}$ of type $[(q^{n-d},d),(1,n-d)]$.}

\bigskip
\noindent{\bf 2.3  Theorem}.  {\it Let $n$, $k$ and $d$ be positive
integers with $d>1$. \break If $n=kd-1$, then there exists a
partition of $V_{n}{(q)}$ of type \break
$[(q^{(k-1)d},d-1),(\frac{q^{(k-1)d}-1}{q^{d}-1},d)]$.}

\bigskip
For partitions of finite vector spaces, we have the following
fundamental problem.

\bigskip
\noindent{\bf 2.4}.  {\it Give necessary and sufficient conditions
on non-negative solutions of the Diophantine equation (1) such that
they correspond to $[(x_{1},n_{1}),(x_{2},n_{2}), \break
......,(x_{k},n_{k})]$-partitions of $V_{n}{(q)}$}.
\bigskip

For small value of $k$ some result is available.

\bigskip
\noindent{\bf 2.5 Proposition}.  {\it The properties  1) and 2) are
necessary and sufficient conditions for the existence of an
$[(x_{1},n_{1}),(x_{2},n_{2}),......,(x_{k},n_{k})]$-partition of
$V_{n}{(q)}$ when $k=1$ or $k=2$ and $ n_{1}+n_{2}=n $}.

\bigskip
{\bf Proof.} \noindent In the case where $k=1$, to every solution of
(1) corresponds an $[(x_{1},n_{1})]$-partition. In fact, the
equation $(q^{n_{1}}-1)x_{1}=q^{n}-1$ admits a solution if and only
if $n_{1}$ divides $n$ and so, by i) of Theorem 2.2, there exists an
$[(x_{1},n_{1})]$-partition.

\noindent Now consider the case where $k=2$ and $n_{1}+n_{2}=n$.
Suppose that $(x_{1},x_{2})$ is a non-negative solution of (1) which
also verifies the necessary condition 2). Since
$n_{2}=n-n_{1}>n_{1}$, we have $n_{1} < \frac{n}{2}$ and so $n_{2} >
\frac{n}{2}$. It follows by 2) that $x_{2} \leq 1$. If $x_{2}=0$,
then from the equation (1) we obtain that $n_{1}$ divides $n$ and
so, again by i) of Theorem 2.2, we have a
$[(\frac{q^{n}-1}{q^{n_{1}}-1},n_{1}),(0,n_{2})]$-partition. In case
$x_{2}=1$, from the equation (1) we get $x_{1}=q^{n-n_{1}}$.
Therefore, $(x_{1},x_{2})=(q^{n-n_{1}},1)$ and so, by ii) of Theorem
2.2, there exists a $[(q^{n-n_{1}},n_{1}),(1,n_{2})]$-partition. So
the proposition is proved.

\bigskip

For $k=2$ and for any $n_{1}$ and $n_{2}$ the question is still an
open problem. Recently, in [ESSSV1] and [ESSSV2] the authors
resolved  it in the case where $n_{1}=2$, $n_{2}=3$ and $q=2$.

\bigskip
\noindent{\bf 2.6  Theorem}.  {\it There exists a partition of
$V_{n}{(2)}$, $n \geq 3$, of type $[(x_{1},2),(x_{2},3)]$ if and
only if $(x_{1},x_{2})$ is a solution of the Diophantine equation

$$ 3x_{1}+7x_{2}=2^{n}-1         \eqno(2) $$
with $x_{1}$ and $x_{2}$ non-negative integers and $x_{1} \neq 1$.}

\bigskip

 In order to show Theorem 2.6, they gave
 the following theorems.

\bigskip
\noindent{\bf 2.7  Theorem}.  {\it Let  $V$ and $V'$ be
$\mathbb{F}_{q}$-vector spaces of finite dimension and
$T=\{n_{1},n_{2},......,n_{k}\}$ a set of positive integers with $
n_{1} < n_{2} < ...... < n_{k}$. If ${\bf P}$ is a $T$-partition of
$V$ with $ n_{k} \leq dim V'$, then there exists a $T$-set of
subspaces ${\bf \overline{P}}$ of $V \oplus V'$ such that $|{\bf
\overline{P}}| = (q^{dim V'}-1)|{\bf P}|$ and $\{ V,V' \} \cup {\bf
\overline{P}} $ is a partition of $V \oplus V'$.}

\bigskip
\noindent{\bf 2.8  Theorem}.  {\it Let ${\bf P}$ be a non-trivial
partition of $V_{n}{(q)}$ of type \break
$[(x_{1},n_{1}),(x_{2},n_{2}),......,(x_{k},n_{k})]$.

i)  If $x_{1} \neq 0$, then $x_{1}\geq 2$.

ii) If  $n_{1}=1$ and $q=2$, then $x_{1}\geq 3$ .}

\bigskip
\noindent{\bf 2.9  Remark}. Note that Theorem 2.7 is a
generalization of Theorem 2.2, part ii). In fact if $d = dim V $,
   $n-d = dim V'$ with $d < n-d $ and ${\bf P} = \{ V \} $, then there exists a partition of
$ V_{n}{(q)} = V \oplus V'$ of type $[(q^{n-d},d),(1,n-d)]$.
Moreover note that Theorem 2.8 gives information about the number of
minimum dimension components  of a partition.

\bigskip
\noindent{\bf 2.10  Remark}. Observe that, if ${\bf P}$ is a
non-trivial $[(x_{1},2),(x_{2},3)]$-partition of $V_{n}{(2)}$ with
$x_{1}$ and $x_{2}$ positive integers, then $x_{1} \geq 3$. In fact,
if $x_{1}=2$, from equation (2), we have $x_{2}=\frac{2^{n}-1-6}{7}$
and so $7$  divides $2^{n}$ which is a contradiction. Of course, if
either $x_{1}$ or $x_{2}$ is equal to zero, we get $x_{1}+x_{2} \geq
3$ since a group can not be the union of two proper and disjoint
subgroups. Therefore, if ${\bf P}$ is a non-trivial
$[(x_{1},2),(x_{2},3)]$-partition of $V_{n}{(2)}$, then the number
of its minimum dimension subspaces is greater than or equal to $3$.

\vspace{1in}

\section{A new necessary condition} \label{S:P*}

\bigskip

\noindent In this section we will prove that, for every prime power
$q$ and for every non-trivial partition ${\bf P}$  of $V_{n}{(q)}$,
the number of minimum dimension subspaces in ${\bf P}$ is always
greater than or equal to $q+t$.

  First we need the following lemmas.

\bigskip
\noindent{\bf 3.1  Lemma}.  {\it Let ${\bf
P}=\{V_{1},V_{2},......,V_{r}\} $  be a non-trivial partition of
$V_{n}{(q)}$ and $t$ and $s$ be positive integers with $t<n$ and
$s<r$. If $dim(V_{i})=t$ for every $i=1,...,s$ and  $dim(V_{j}) \geq
t+1$ for every $j=s+1,...,r$, then $s \geq \alpha q $ for some
positive integer $ \alpha $.}

\bigskip
{\bf Proof.} We can suppose that $V_{n}{(q)}=\mathbb{F}_{q}^ {n}$
and $V_{1}=\mathbb{F}_{q} ^{t} \times \{ 0 \}^{n-t}$ after we choose
an ordered basis of $V_{n}{(q)}$ which extends a fixed ordered basis
of $V_{1}$. Let $W=V_{1}^ {\perp}$ be the dual subspace of $V_{1}$
with respect to  the canonical inner product on $\mathbb{F}_{q}^
{n}$. Of course we have $V_{1} \cap W=\{ 0 \}$ and $V_{n}{(q)}=V_{1}
\oplus W$. Moreover, for every $j=s+1,...,r$, we have that
$dim(V_{j}+W)=dim(V_{j})+dim(W)-dim(V_{j} \cap
W)=n-t+dim(V_{j})-dim(V_{j} \cap W)$. But $dim(V_{j})>dim(V_{1})=t$,
then  $z:=dim(V_{j})-t \geq 1$. Thus we have
$dim(V_{j}+W)=n+z-dim(V_{j} \cap W)\leq n$. It follows $dim(V_{j}
\cap W) \geq z \geq 1$ and so we get $V_{j} \cap W \neq \{0\}$.

\noindent Of course it is possible that there are some other
subspaces of dimension $t$, different from $V_{1}$, which are not
disjoint from $W$. So, we can suppose that there exists an integer
$s'$ with $1 \leq s' \leq s$ and such that $V_{i} \cap W = \{0\}$
for every $i=1,...,s'$, whereas $V_{i} \cap W \neq \{0\}$ for every
$i>s'$ and $i \leq s$. Consider the partition ${\bf P'}$ induced by
${\bf P}$ on $W$, that is ${\bf P'}=\{ V_{i} \cap W \hspace{0.05in}
| \hspace{0.05in} V_{i} \in {\bf P} \hspace{0.05in} and
\hspace{0.05in} V_{i}\cap W \neq \{0\} \}$. Clearly we have

 $$ \sum_{j=s'+1}^{r}(q^{m_{j}}-1)= q^{n-t}-1, \eqno
(4)$$
 where $m_{j}=dim(V_{j} \cap W) \geq 1$ for every $j=s'+1,...,r$.
Further, considering the partition ${\bf P}$, we obtain $
\sum_{j=s'+1}^{r}(q^{n_{j}}-1)= q^{n}-1-
\sum_{j=1}^{s'}(q^{n_{j}}-1)=q^{n}-1-s'(q^{t}-1)$, from which

$$ \sum_{j=s'+1}^{r}(q^{n_{j}}-1)= q^{n}-1-s'q^{t}+s'.
\eqno (5)
$$
\noindent Now, subtracting  (4) from (5), we get

$$ \sum_{j=s'+1}^{r}(q^{n_{j}}-q^{m_{j}})= q^{n}-s'q^{t}-q^{n-t}+s'. $$

\noindent Hence we obtain that $q$ divides $s'$. But $s' \neq 0$
being  $V_{1} \cap W=\{ 0 \}$. Therefore, for some positive integer
$\alpha$, we obtain $s \geq s' = \alpha q$ and the proof is
complete.

\bigskip
\noindent{\bf 3.2 Lemma}.  {\it Let ${\bf P}$ be a non-trivial
partition of $V_{n}{(q)}$ of type \break
$[(x_{1},n_{1}),(x_{2},n_{2}),......,(x_{k},n_{k})]$. Suppose
$i_{0}$ be a positive integer such that $ x_{i_{0}} \neq 0 $ and $
x_{i}=0 $ for each $ i=1,2,...,i_{0}-1 $, then $ x_{i_{0}} \geq
\alpha q+1 $ where $\alpha$ is a positive integer.}

\bigskip
{\bf Proof.} We use the same notations as in the previous lemma. If
$i_{0}=k$, then we have that $ x_{i_{0}}=x_{k}= \frac
{q^{n}-1}{q^{n_{k}}-1}$ because of condition 1). So $n_{k}$ divides
$n$. Moreover, being ${\bf P}$ a non-trivial partition, we have $
n_{k} < n $. It follows that
$x_{i_{0}}=q^{n-n_{k}}+q^{n-2n_{k}}+.....+q^{n-(\frac
{n}{n_{k}}-1)n_{k}}+1= \alpha q^{n_{k}}+1 \geq \alpha q+1$ where
$\alpha =q^{n-2n_{k}}+q^{n-3n_{k}}+.....+q^{n_{k}}+1$.

\noindent Now suppose $i_{0} < k$ and set $t=n_{i_{0}}$ and
$s=x_{i_{0}}$ as in the previous lemma. If ${\bf
P}=\{V_{1},V_{2},......,V_{r}\}$, and $V_{1},V_{2},......,V_{s}$ are
the $s$ components of dimension $t$, then there exist at least two
distinct components which have the same  minimum dimension $t$
since, by Lemma 3.1, $s \geq q \geq 2$. Therefore, we certainly have
that $V_{1}$ is distinct from $V_{s}$. So $dim(V_{1}+V_{s})=2t \leq
n$ because of $V_{1}$ and $V_{s}$ belong to ${\bf P}$ and
$dim(V_{1})=dim(V_{s})=t$. Let $\{ v_{1},v_{2},.....,v_{t} \}$ be an
ordered basis of $V_{1}$ and $\{ v'_{1},v'_{2},.....,v'_{t} \}$ be
an ordered basis of $V_{s}$. Then the vectors $\{
v_{1},v_{2},.....,v_{t},v'_{1},v'_{2},.....,v'_{t} \}$ are a basis
of $V_{1}+V_{s}$. So they are linearly independent and we can
consider a basis of $V_{n}{(q)}$ which contains them. In relation to
this new basis we can identify $V_{n}{(q)}$ with
$\mathbb{F}_{q}^{n}$, $V_{1}$ with $\mathbb{F}_{q}^{t} \times
\{0\}^{n-t}$ and $V_{s}$ with $\{0\}^{t} \times \mathbb{F}_{q}^{t}
\times \{0\}^{n-2t}$. Let $W$ be the dual space of $V_{1}$, that is
$W=\{0\}^{t} \times \mathbb{F}_{q}^{n-t}$. It follows that $V_{s}
\subseteq W$ and so $V_{s} \cap W \neq \{0\}$. Therefore, the number
$s'$ of the $t$-dimensional subspaces of $V_{n}{(q)}$ which  are
disjoint from $W$ is smaller than $s$. But, as in the proof of Lemma
3.1, we have $s' = \alpha q$ for some  positive integer $\alpha$. So
$s
> s' = \alpha q$, that is to say $x_{i_{0}}=s \geq \alpha q+1$, and
the lemma is shown.

\bigskip
\noindent{\bf 3.3 Theorem}. {\it In any non-trivial partition of
$V_{n}{(q)}$, the number of subspaces of minimum dimension $t$ is
greater than or equal to $ \alpha q+t $ for some positive integer
$\alpha$.}

\bigskip

{\bf Proof.} We proceed by induction on $t$. For  $t = 1$ the
theorem is true by the above Lemma 3.2. Now, let $t \geq 2$, ${\bf
P}$ be a partition of $V_{n}{(q)}$ and $S$ be the subset of ${\bf
P}$ of all the components of minimum dimension $t$.  Consider an
hyperplane $W$ of $V_{n}{(q)}$ which contains at least a component
of $S$ and it does not contain all the contains of $S$. A such
hyperplane there exists since $s = |S| \geq \alpha q+1 \geq q+1 \geq
3$ by the above lemma. Being $t
> 1$ the partition ${\bf P}_{W}$, which is induced from ${\bf P}$ on
$W$, has components of minimum dimension $t-1$. Let $S'$ be the set
of such components of ${\bf P}_{W}$. So, by induction, their number
is $s' = |S'| = \alpha q+t-1$. But if $V'_{i} \in S'$, then $V'_{i}
= V_{i} \cap W$ where $V_{i} \in S$. Therefore, $s \geq s' = \alpha
q+t-1$. Moreover, by construction, $W$ contains at least one
component of $S$. Thus we get that
 $s \geq s'+1 \geq
(\alpha q+t-1)+1$. So $s \geq \alpha q+t$ and the proof is complete.

\bigskip

Finally we can state the next corollary which clearly follows from
the above theorem.

\bigskip
\noindent{\bf 3.4 Corollary}. {\it Let $V_{n}{(q)}$ be a vector
space which admits a non trivial partition ${\bf P}$. Then the
number of components of ${\bf P}$ of minimum dimension $t$ is
greater than or equal to $ q+t $ .}

\bigskip
 Now we observe that, if ${\bf P}=\{V_{1},V_{2},......,V_{r}\} $ is a
non-trivial partition of $V_{n}{(q)}$ whose components are all of
the same dimension $t$, then  $r$ is equal to the number $s$ of
minimum dimension components of ${\bf P}$ and $t$ divides $n$ by
Proposition 2.5. So we have that $r = s = \frac{q^{n}-1}{q^{t}-1}
\geq q^{t}+1$ (Note that, $s$ may be much  greater than $q+t$ if $t
\neq 1$). More generally, we have the following proposition.

\bigskip
\noindent{\bf 3.5 Proposition}.  {\it Let  ${\bf P}$  be a
non-trivial partition of $V_{n}{(q)}$ which have $r$ components. If
$t$ is the minimum dimension of the components of ${\bf P}$, then
$q^{t}+1 \leq r \leq \lfloor \frac{q^{n}-1}{q^{t}-1} \rfloor $.

\noindent (Here $ \lfloor x \rfloor $ denotes the integer part of
the real number $x$) .}

\bigskip
{\bf Proof.} Suppose ${\bf P}=\{V_{1},V_{2},......,V_{r}\}$ be the
partition of $V_{n}{(q)}$. Then we have $
\sum_{i=1}^{r}(q^{n_{i}}-1)= q^{n}-1$ if $n_{i}$ denotes the
dimension of $V_{i}$ for every $1 \leq i \leq r$. Hence $r-1 =
\sum_{i=1}^{r} q^{n_{i}}-q^{n}= q^{t}(\sum_{i=1}^{r}
q^{n_{i}-t}-q^{n-t})$ and we obtain that  $r = \alpha q^{t}+1$ with
$\alpha \geq 1$ being ${\bf P}$ a non-trivial partition. It follows
that $r \geq q^{t}+1$. Now, let $ \{ V_{1},V_{2},......,V_{s} \}$ be
the components of ${\bf P}$ of minimum dimension $t$ and suppose $s
< r$. We have that $(V \backslash \{ 0 \}) \backslash (
\bigcup_{i=1}^{s}(V_{i} \backslash \{ 0 \} ))=
\bigcup_{i=s+1}^{r}(V_{i} \backslash \{ 0 \} )$. So we obtain $|
\bigcup_{i=s+1}^{r}(V_{i} \backslash \{ 0 \} ) | = | V \backslash \{
0 \} | - | \bigcup_{i=1}^{s}(V_{i} \backslash \{ 0 \} ) |$. But
$(r-s)(q^{t}-1) \ < | \bigcup_{i=s+1}^{r}(V_{i} \backslash \{ 0 \} )
|$ since $| V_{i} | > q^{t}$ for every $s+1 \leq i \leq r$. It
follows that $(r-s)(q^{t}-1) < | V \backslash \{ 0 \} | - |
\bigcup_{i=1}^{s}(V_{i} \backslash \{ 0 \} ) | =
(q^{n}-1)-s(q^{t}-1)$ and so $(r-s)(q^{t}-1) < q^{n}-1-sq^{t}-s$
from which  we get $r < \frac{q^{n}-1}{q^{t}-1}$. If $s = r$, then
the components of ${\bf P}$ have all the same dimension $t$ and, as
observed before, $r = \frac{q^{n}-1}{q^{t}-1}$. Therefore, in any
case, $r \leq \frac{q^{n}-1}{q^{t}-1}$ and the proposition is shown.

\bigskip
\noindent{\bf 3.6 Remark}. Proposition 3.5 and the examples of
partitions which are known to us, drive us to think that Corollary
3.4 can be substantially improved. In fact, we conjecture that the
number of components of minimum dimension $t$ of a non-trivial
partition  of $V_{n}{(q)}$ is greater or equal to $q^{t}+1$.

 \vspace{1in}

\section{Existence results on $T$-partitions} \label{S:P*}

\bigskip

In this section we give some extensions of Theorem 2.1. To begin, we
can drop the hypothesis "$ n_{1} \geq 2 $" in Theorem 2.1. In fact,
we have the following proposition.

\bigskip
\noindent{\bf 4.1 Proposition}.  {\it Let $T = \{
n_{1},n_{2},......,n_{k} \}$ be a set of positive integers such that
$ n_{1}< n_{2}<......< n_{k}$. Then there exists a T-partition of
$V_{2n_{k}}{(q)}$.}

\bigskip
{\bf Proof.} By Theorem 2.1 we can suppose that $ n_{1} = 1 $. If $k
= 1$ the proposition follows  by $i)$ of Theorem 2.2.  So let $k
\geq 2$ and consider the subset $T'= \{ n_{2},n_{3}......,n_{k} \} $
of $T$. Again by Theorem 2.1, there exists a $T'$-partition ${\bf
P'}$ of $V_{2n_{k}}{(q)}$ because $n_{2}
 > n_{1} = 1$. But, by Corollary 3.4, there exist at least $q+n_{2} \geq 4$
 components of ${\bf P'}$ of minimum dimension $n_{2}$. Let $V'$ be
 such a component of dimension $n_{2}$ and
 consider the partition  ${\bf P''}$ of $V'$ whose components are
 all its  subspaces of dimension 1.
  Now, ${\bf P} = ({\bf P'} \setminus \{ V' \}) \cup {\bf P''}$
  is a $T$-partition of $V_{2n_{k}}{(q)}$ since $|{\bf P''}| \geq q+1 \geq 1$
  and there are some other components (at least 3) of ${\bf P'} \subset {\bf P}$
  of dimension $n_{2}$. This complete the proof.

\bigskip

\noindent{\bf 4.2 Theorem}.  {\it Let  $T = \{
n_{1},n_{2},....,n_{k} \}$ be a set of positive integers such that $
n_{1} < n_{2} < .... < n_{k} $   and consider the vector space $
V_{n}{(q)}$ with  $n \geq 2n_{k}$. If  $gcd(n,2n_{k})$ admits some
divisor into $T$,
 then there exists a
$T$-partition of $ V_{n}{(q)}$.}

\bigskip
{\bf Proof.} By Proposition 4.1 we can suppose that $n > 2n_{k}$.
 Consider a subspace $V$ of $V_{n}{(q)}$ of dimension $2n_{k}$ and
 such that $V \cap V^{\perp} = \{ 0 \}$ where $V^{\perp}$ is the
 dual space of $V$. If $n_{i_{0}} \in T$ is a divisor of $gcd(n,2n_{k})$, then
$n_{i_{0}}$ is a divisor of $n - 2n_{k} = dim V^{\perp} $. So, by
Theorem 2.2, there exists a partition ${\bf P'}$ of $V^{\perp}$
whose components have all the same dimension $n_{i_{0}}$. Since
$n_{i_{0}}$ divides $2n_{k}$, for the same reason we get that there
exists a partition ${\bf P}$ of $V$ whose components have all the
same dimension $n_{i_{0}}$, that is, ${\bf P}$ is a
$\bar{T}$-partition of $V$ where  $\bar{T} = \{ n_{i_{0}} \}$ and
$n_{i_{0}} \leq dim V^{\perp} = n-2n_{k} $ being $n_{i_{0}}$ a
divisor of $n-2n_{k}$. It follows, by Theorem 2.7, that there exists
a $\bar{T}$-set
 of $n_{i_{0}}$-dimensional subspaces ${\bf P''}$ of $V \oplus V^{\perp} =
V_{n}{(q)}$ such that $  \{ V,V^{\perp} \} \cup {\bf P''}$ is a
partition of $ V_{n}{(q)}$. Now, by the above proposition, let
$\bar{{\bf P}}$ be a $T$-partition of $V$. Then $\bar{{\bf P}} \cup
{\bf P'} \cup {\bf P''}$ is a $T$-partition of $V_{n}{(q)}$ and the
proof is complete.

\bigskip

\noindent Note that the above theorem is Proposition 4.1 for
$n=2n_{k}$. So it is a generalization of Beutelspacher-Heden's
Theorem 2.1.

\bigskip

\noindent{\bf 4.3 Lemma}.  {\it Let $T$ be as in Theorem 4.2 and $
V_{n}{(q)}$ be a vector space over $\mathbb{F}_{q}$ of dimension $n
\geq 3n_{k}$. If there exists a subset $T'$ of $T$ such that
$V_{n-2n_{k} }{(q)}$ has a $T'$-partition, then there exists a
$T$-partition of $ V_{n}{(q)}$ }.

\bigskip

{\bf Proof.} Let $V$ be a subspace of $ V_{n}{(q)}$ of dimension
$2n_{k}$ and such that $V \cap V^{\perp} = \{ 0 \}$. Since $dim
V^{\perp} = n-2n_{k} $, then $V^{\perp}$ is (isomorphic to)
$V_{n-2n_{k} }{(q)}$ and so $V^{\perp}$ admits a $T'$-partition
${\bf P'}$. By Proposition 4.1, we can consider a $T$-partition
${\bf P}$ of $V$. Moreover, by hypothesis, $n-2n_{k} \geq n_{k}$ and
so $dim V^{\perp} \geq n_{k}$. Therefore there exists a $T$-set of
subspaces ${\bf P''}$ of $V \oplus V^{\perp}$ such that $ \{
V,V^{\perp} \} \cup {\bf P''}$ is a partition of $V \oplus
V^{\perp}$. Now we get that $ {\bf P} \cup {\bf P'} \cup {\bf P''}$
is a $T$-partition of $V \oplus V^{\perp} = V_{n}{(q)}$ and the
lemma is shown.

\bigskip

For $n$ greater than or equal to $3n_{k}$ or for $n$ smaller than
$2n_{k}$, we have the next theorem.

\bigskip
 \noindent{\bf 4.4 Theorem}.  {\it Let $n$ be a positive integer and $T = \{
n_{1},n_{2},....,n_{k} \}$ be a set of positive integers such that $
n_{1} < n_{2} < .... < n_{k} < n $. Then  $ V_{n}{(q)}$ admits a
$T$-partition if one of the following hypothesis is satisfied:

a)  $ 3n_{k} \leq n$ and $n-2n_{k}$ admits some divisor into $T$;

b) $2n_{k} > n = n_{k}+n_{k-1}$ and $n_{1}=1$;

c) $2n_{k} > n \geq n_{k}+2n_{k-1}$ and $gcd(n,2n_{k-1})$ has some
divisor into $T$.}

\bigskip
{\bf Proof.} a) By the above lemma it is enough to note that there
exists $T' \subseteq T$ such that $V_{n-2n_{k} }{(q)}$ has a
$T'$-partition. In fact, if  $ n_{i_{0}} \in T$ is a divisor of
$n-2n_{k}$, then $V_{n-2n_{k} }{(q)}$ admits a partition whose
components have all the same dimension $n_{i_{0}}$. That is to say,
$V_{n-2n_{k} }{(q)}$ has a $T'$-partition  if we set $T' = \{
n_{i_{0}} \}$.

b) Note that, since $n_{k} > \frac{n}{2}$, if $ V_{n}{(q)}$ admits a
$T$-partition then exists exactly one subspace of dimension $n_{k}$
because of the necessary condition 2). Again we note that if $d$ is
a positive integer smaller than the dimension of a vector space $U$,
then always there exists a $\bar{T}$-partition of $U$ where $\bar{T}
= \{ 1, d \}$; in fact, the components of a such $\bar{T}$-partition
may be a fixed $d$-dimensional subspace $W$ of $U$ and the
$1$-dimensional subspaces which are not in $W$.

\noindent Now we can prove b). By ii) of Theorem 2.2, let ${\bf P'}$
be a $T'$-partition of $ V_{n}{(q)}$ of type
$[(q^{n_{k}},n-n_{k}),(1,n_{k})]$. Of course here $T' = \{
n-n_{k},n_{k} \}$. Since $q^{n_{k}} \geq n_{k} \geq k > k-2$,  we
can choose $k-2$ distinct subspaces $V_{1},V_{2},...,V_{k-2}$ of
dimension $n-n_{k} = n_{k-1}$ of the $T'$-partition ${\bf P'}$. For
every $1 \leq i \leq k-2$, being $n_{i} < n_{k-1} = n-n_{k} = dim
V_{i}$, it is possible to consider
  a $T'_{i}$-partition ${\bf P'_{i}}$ of $V_{i}$ where $T'_{i} = \{ 1, n_{i}
\}$. Let ${\bf P''} = \{ V_{i} \in {\bf P'}
 \hspace{0.05in} | \hspace{0.05in} k-1 \leq i \leq q^{n_{k}} \hspace{0.05in} and \hspace{0.05in} dim V_{i} = n_{k-1} \}$.
 Since $q^{n_{k}} > k-2$, we get that ${\bf P''} \neq \emptyset$.
Now, if $V$ is the component of dimension $n_{k}$ in ${\bf P'}$,
then $\{ V \}  \cup (\cup_{i=1}^{k-2}{\bf P'_{i}}) \cup {\bf P''}$
 is a $T$-partition of $ V_{n}{(q)}$.

c) Let $V$ be a subspace of $ V_{n}{(q)}$ of dimension $n-n_{k}$ and
such that $V \cap V^{\perp} = \{ 0 \}$. The hypothesis $n > n_{k} >
\frac{n}{2} $ implies that $n_{k}$ is not a divisor of
$gcd(n,2n_{k-1})$ and so $gcd(n,2n_{k-1})$ admits a divisor in $T' =
\{ n_{1},n_{2},....,n_{k-1} \}$. Now  since $n-n_{k} \geq 2n_{k-1}$,
by Theorem 4.2, we have that there exists a $T'$-partition ${\bf
P'}$ of $V$.  The subspace $V^{\perp}$ has dimension $n_{k}$ and so
$n_{k-1} < dim V^{\perp}$. Therefore, by Theorem 2.7, we obtain that
there exists a $T'$-set of subspaces ${\bf P''}$ such that $\{ V,
V^{\perp} \} \cup {\bf P''} $ is a partition of $V_{n}{(q)}$. Now we
get that ${\bf P'} \cup {\bf P''} \cup \{ V^{\perp} \} $ is a
$T$-partition of $V_{n}{(q)}$. So the theorem is completely shown.

\bigskip

\noindent{\bf 4.5 Remark}. Of course the above results do not give a
complete answer to the problem to give   sufficient conditions for
the existence of a $T$-partition of  $ V_{n}{(q)}$. For example, it
is not known if there are $T$-partitions of $ V_{n}{(q)}$ when $
3n_{k}
> n > 2n_{k}$ and no integer belonging to $T$ divides $gcd(n,2n_{k})$.

More generally the problem to give necessary and sufficient
conditions about a set of positive integers $T$ to have a
$T$-partition of $ V_{n}{(q)}$ is an open problem. Of course, it is
in correlation with the analogous problem on \break
$[(x_{1},n_{1}),(x_{2},n_{2}),......,(x_{k},n_{k})]$-partitions.
 But, now the problem is to give conditions on the elements of $T$
 such that there exists some \break
 $[(x_{1},n_{1}),(x_{2},n_{2}),......,(x_{k},n_{k})]$-partition
 where $x_{1},x_{2},...,x_{k}$ are positive integers. Note that, as it follows from the results of this section, it
 is not necessary to know the positive  integers $x_{1},x_{2},...,x_{k}$ to have the
 existence of a $T$-partition.

\vspace{1,2in}

\noindent {\Large {\bf References}}

\bigskip
\bigskip

\noindent [Be] Beutelspacher A. {\it Partitions of finite vector
spaces: an application of the Frobenius number in geometry.} Arch.
Math. 31 (1978), 202-208.

\bigskip
\noindent [Bu] BU T. {\it Partitions of a vector space.} Discrete
Math. 31 (1980), 79-83.

\bigskip
\noindent [DR] Danziger P. and Rodney P. {\it Uniformly resolvable
designs } in: The CRC handbook of combinatorial designs. Edited by
Charles J. Colbourn and Jeffrey H. Dinitz. CRC Press Series on
Discrete Mathematics and its applications. CRC Press, Boca Raton,
FL, (1996), 490-492.

\bigskip
\noindent [ESSSV1]  El-Zanati S.I., Seelinger G.F., Sissokho P.A.,
Spence L.E. and Vanden Eynden C. {\it On partitions of finite vector
spaces of small dimension over GF(2).} Submitted.

\bigskip
\noindent [ESSSV2]  El-Zanati S.I., Seelinger G.F., Sissokho P.A.,
Spence L.E. and Vanden Eynden C. {\it Partitions of finite vector
spaces into subspaces.} J. Comb. Des. 16,No.4 (2008) 329-341.

\bigskip
\noindent [He] Heden O.{\it On partitions of finite vector spaces of
small dimensions.} Arch. Math. 43(1984), 507-509.

\bigskip
\noindent [HS]  Herzog M. and Schonheim J. {\it Linear and nonlinear
single error-correcting perfect mixed codes.} Informat. and Control
18 (1971), 364-368.

\bigskip
\noindent [Li]  Lindstrom B. {\it Group partitions and mixed perfect
codes.} Canad. Math. Bull. 18 (1975), 57-60.

\bigskip
\noindent [Sc1] Schulz R.-H.  {\it On existence of Generalized
triads related to transversal designs.} Ars Comb. 25 B (1988),
203-209.

\bigskip
\noindent [Sc2] Schulz R.-H. {\it On translation transversal designs
with $\lambda > 1$.} Arch. Math. 49,(1987) 97-102.

\bigskip
\noindent [Sp] Spera A.G. {\it $(s,k,\lambda)$-partitions of a
vector space.} Discrete Math. 89 (1991), 213-217.

\end{document}